\documentclass{amsart}

\usepackage{amsmath,amsfonts,amsthm,xypic,color}

\newcommand{\ACF}{\mathsf{ACF}}
\newcommand{\FO}{\mathsf{FO}}
\newcommand{\Th}{\mathrm{Th}}
\DeclareMathOperator{\Spec}{Spec}

\newcommand{\mbf}{\mathbf}

\newcommand{\findaref}{}

\newtheorem{thm}{Theorem}
\newtheorem{cor}[thm]{Corollary}
\newtheorem{lem}[thm]{Lemma}
\newtheorem{prop}[thm]{Proposition}
\theoremstyle{definition}

\newtheorem{ex}[thm]{Example}
\newtheorem{rmk}[thm]{Remark}

\title[The first-order theory of geometric points of schemes]{The first-order theory of geometric points of schemes: Chevalley's theorem and quantifier elimination}
\author{L.\ Alexander Betts}
\address{Merton College, Oxford OX1 4JD, United Kingdom}
\email{alexander.betts@maths.ox.ac.uk}

\begin{document}

\maketitle

\vspace{-1cm}

\begin{abstract}
Chevalley's theorem on the images of morphisms of schemes and the principle of quantifier elimination for the theory of algebraically closed fields are widely understood to be two perspectives on the same theorem. In this paper, we demonstrate that both results can easily be proven simultaneously, using a mixture of geometric and logical techniques. In doing so, we give logical meaning to geometric points of schemes and to finitely presented morphisms thereof, in a manner reminiscent of Spencer Breiner's logical schemes.
\end{abstract}

\section*{Introduction}

In his PhD thesis \cite{breiner}, Breiner develops the notion of a logical scheme. To a good approximation, this is a geometric object which is locally a first-order theory, and its points correspond to models of these theories.

We can realise a similar idea concretely in the case of (classical) schemes. To each ring $R$ we can associate a first-order theory $\ACF_R$, so that we have a canonical bijection\[\{\text{Geometric points of $\Spec(R)$}\}\longleftrightarrow\{\text{Models of $\ACF_R$}\}\]This correspondence induces a bijection between the $0$th Stone space of $\ACF_R$ and the underlying topological space of $\Spec(R)$, which is a homeomorphism when $|\Spec(R)|$ is endowed with the constructible topology.


Using this correspondence we can easily understand the images of (locally) finitely presented morphisms of schemes:  they are just (locally) described by a sentence in the language of $\ACF_R$. Pursuing this idea naturally leads to a straightforward proof of Chevalley's theorem (on the scheme-theoretic side) and quantifier elimination for the theories $\ACF_R$ (on the logical side), and the aim of this paper is to explain these basic applications of the theory.

\begin{rmk}
Although our treatment involves heavy use of the underlying topological space of a scheme, for the purely algebraic applications of this theory it is easy to dispense with this and just work with the functor-of-points instead.
\end{rmk}

\section*{Geometric points of schemes as models of theories}

The correspondence for geometric points is very simple. We form the language of $\ACF_R$ by taking the language of rings and adjoining constants $c_a$ for each element $a\in R$ (we will frequently write $c_a$ just as $a$, for ease of notation). The theory $\ACF_R$ is then formed by taking the theory of algebraically closed fields (in the language of rings) and adjoining axioms $c_1=1$ and, for each $a,b\in R$, $c_{a+b}=c_a+c_b$ and $c_{ab}=c_ac_b$.

Models of $\ACF_R$ are then algebraically closed fields $k$, together with named elements $c_a^k$ for each $a\in R$, so that the function $R\rightarrow k$ sending $a$ to $c_a^k$ is a homomorphism of rings. In other words, models of $\ACF_R$ can be naturally thought of as ring homomorphisms $x:R\rightarrow k$ from $R$ into an algebraically closed field (or more accurately, such a homomorphism makes $k$ into a model of $\ACF_R$ in a natural way -- we will be happy to conflate these concepts).

Since $\Spec$ is a contravariant embedding, such homomorphisms are in bijection with scheme morphisms $\Spec(k)\rightarrow\Spec(R)$, which are precisely geometric points. Thus the correspondence between geometric points of $\Spec(R)$ and models of $\ACF_R$ is just transposition across the $\Spec$/global sections adjunction.

\subsection*{The image formula}

With this correspondence between geometric points and models of $\ACF_R$ in mind, we make the crucial observation that the image of a finitely presented morphism of affine schemes is described by a particular first-order formula. It follows that the image will be constructible iff this formula can be rewritten in an equivalent, quantifier-free form. Of course, elimination of quantifiers for the theories $\ACF_R$ makes this immediate, but since we are aiming to prove this, we use more elementary methods.

\begin{prop}\label{prop:image_formula}
Let $S\cong R[\mbf t]/(\mbf f)$ be a finitely presented $R$-algebra, and $\Spec(S)\rightarrow\Spec(R)$ the structure map. Then a geometric point $x:\Spec(k)\rightarrow\Spec(R)$ lifts to $\Spec(S)$ iff, viewed as a model of $\ACF_R$, it satisfies the sentence\[(\exists\mbf y)\left(\bigwedge_j(f_j(\mbf y)=0)\right)\]

We will refer to this formula as the \emph{image formula} associated to this presentation of $S$ as an $R$-algebra.
\begin{proof}
The lifts of $x$ to $\Spec(S)$ correspond with homomorphisms $R[\mbf t]/(\mbf f)\rightarrow k$ of $R$ algebras, and hence with tuples $\mbf y\in k^n$ such that all $f_j(\mbf y)=0$. Thus $x$ lifts to $\Spec(S)$ iff such a tuple exists, which says precisely that $x:R\rightarrow k$ satisfies the image formula.
\end{proof}
\end{prop}

The preceding proposition gives logical meaning to the collection of geometric points in $\Spec(R)$ lifting to $\Spec(S)$, but we should check that this agrees with the usual topological notion of the image of $\Spec(S)\rightarrow\Spec(R)$. This is provided (in a strong sense) by the following easy proposition, which is also the sole algebraic input into the theory.

\begin{prop}\label{prop:geom_image=top_image}
Let $Y\rightarrow X$ be a locally finitely presented morphism of schemes, and $x:\Spec(k)\rightarrow X$ a geometric point with underlying topological point $P\in|X|$. Then $x$ lifts to $Y$ iff $P$ lies in the image of $Y\rightarrow X$. In particular, whether $x$ lifts to $Y$ depends only on its image in $X$.
\begin{proof}
Let $\kappa_P$ be the residue field at $P$, so that $x$ factors (uniquely) through $\Spec(\kappa_P)\rightarrow X$. We know that $P$ lies in the image of $Y\rightarrow X$ iff the fibre $\Spec(\kappa_P)\times_XY$ is not the empty scheme. Since $\kappa_P\rightarrow k$ is faithfully flat, this occurs iff $\Spec(k)\times_XY$ is not the empty scheme. Yet this is locally finitely presented over $\Spec(k)$, so the Nullstellensatz shows that $\Spec(k)\times_XY$ is non-empty iff its structure morphism to $\Spec(k)$ has a section. By the universal property of fibre products, such sections are in bijective correspondence with lifts of $x$ to $Y$, completing the proof.
\end{proof}
\end{prop}

\begin{ex}
If $R$ is any ring and $b\in R$ an element, then the (finitely presented) morphism $\Spec(R_b)\rightarrow\Spec(R)$ is just the inclusion of a standard Zariski-open affine. Under $\ACF_R$, its image formula is equivalent to the formula $(b\neq0)$, so that the image formula describing $D(b)$ is just the statement that $b$ does not vanish.
\end{ex}

\section*{Chevalley's theorem}

We are now in a position to prove Chevalley's theorem on the images of finitely presented scheme morphisms. The argument is essentially a two-pass application of the compactness theorem from first-order logic, but the niceties can be cut out by recognising the argument as exactly the same one that proves the transfer lemma\findaref{} from model theory. For clarity, we recall this here before applying it to prove the affine case of Chevalley's theorem.

\begin{lem}[Transfer lemma]\label{lem:transfer_lemma}
Let $T_0\subseteq\FO_0$ be some base theory (over some signature) and $F\subseteq\FO_0$ a fragment closed under $\vee$ and $\wedge$. Suppose that a sentence $\Phi$ has the following property: whenever $A$ and $B$ are two $T_0$-models such that $A\models\Phi$ and $F\cap\Th(A)\subseteq\Th(B)$, then $B\models\Phi$ also. Then $\Phi$ is equivalent under $T_0$ to a sentence in $F$.
\end{lem}

\begin{thm}[Chevalley's theorem (affine case)]\label{thm:affine_Chevalley}
Let $S$ be a finitely presented $R$-algebra. Then the image of the map $|\Spec(S)|\rightarrow|\Spec(R)|$ is constructible. Equivalently, the image formula $\Phi$ for $\Spec(S)\rightarrow\Spec(R)$ is equivalent under $\ACF_R$ to a quantifier-free sentence.
\begin{proof}
Translating proposition \ref{prop:geom_image=top_image} into the language of model theory, we see that whether a model $x:R\rightarrow k$ of $\ACF_R$ satisfies the image formula $\Phi$ depends only on $\ker(x)\unlhd R$. In particular, whether a model satisfies $\Phi$ depends only on which quantifier-free sentences it satisfies. Thus $\Phi$ and the quantifier-free fragment of $\FO_0$ satisfy the conditions of the transfer lemma, where we take $T_0=\ACF_R$. As a consequence, we see that $\Phi$ is equivalent under $\ACF_R$ to a quantifier-free sentence, as desired.

To pass back to constructibility is now easy. $\Phi$ is equivalent to a sentence of the form $\bigvee_i\left((b_i\neq0)\wedge\bigwedge_j(a_{ij}=0)\right)$. Translating into the language of scheme theory, a geometric point of $\Spec(R)$ lies in the image of $\Spec(S)\rightarrow\Spec(R)$ iff its underlying topological point lies in $\bigcup_i\left(D(b_i)\cap\bigcap_j V(a_{ij})\right)$, which is a constructible set.
\end{proof}
\end{thm}


\begin{rmk}
Exactly the same method can be used to prove that finitely presented flat morphisms $f:\Spec(S)\rightarrow\Spec(R)$ of affine schemes have Zariski-open image. A standard result \cite[Theorem 5.D]{matsumura} in commutative algebra tells us that the topological image of $f$ is closed under generisation. In other words, given geometric points $x:R\rightarrow k$ and $y:R\rightarrow l$ such that $\ker(y)\leq\ker(x)$, if $x$ lies in the image of $f$ then so does $y$. Thus we may apply the transfer lemma immediately, with respect to the fragment of $\FO_0$ consisting of the quantifier-free \emph{negative} sentences, to deduce that the image formula $\Phi$ of $f$ is equivalent under $\ACF_R$ to a quantifier-free negative sentence. In other words, the image is Zariski-open (indeed is a finite union of basic opens).
\end{rmk}

\begin{cor}[Chevalley's theorem]
Let $f:Y\rightarrow X$ be a locally finitely presented, quasicompact morphism of schemes. Then the image of any locally constructible set in $Y$ is locally constructible in $X$.
\begin{proof}
We can immediately reduce to the case when $Y=\Spec(S)$ and $X=\Spec(R)$ are affine, so that the morphism is given by a finitely presented ring homomorphism $R\rightarrow S$. Any constructible set $Z$ in $\Spec(S)$ is the image of a finitely presented morphism $\Spec(T)\rightarrow\Spec(S)$ and thus $f(Z)$ is the image of the composite $\Spec(T)\rightarrow\Spec(R)$. It is hence constructible by theorem \ref{thm:affine_Chevalley}.
\end{proof}
\end{cor}

\begin{cor}[Quantifier elimination for $\ACF_R$]\label{cor:QE_for_ACF_R}
The theory $\ACF_R$ has quantifier elimination. Every first-order formula in the language of $\ACF_R$ is equivalent to a quantifier-free formula.
\begin{proof}
Proceeding by structural induction, we just need to show that existential formulae are equivalent to quantifier-free formulae. It is easy to see that every quantifier-free formula $\Psi$ in the language of $\ACF_R$ is equivalent to one in the form\[\bigvee_i\left((b_i\neq0)\wedge\bigwedge_j(a_{ij}=0)\right)\]for some $R$-polynomials $b_i,a_{ij}$ in the free variables of $\Psi$.

If the formula $(\exists\mbf x)\Psi$ has no free variables, it is equivalent to the image formula of some finitely presented morphism of affine schemes (associated to the $R$-algebra $\prod_iR[\mbf t]_{b_i}/(a_{ij})_j$), so that we are done by theorem \ref{thm:affine_Chevalley}. In general, if the formula $(\exists\mbf x)\Psi$ has free variables $\mbf y$, we can view it as a sentence in the language of $\ACF_{R[\mbf y]}$, whose models can be thought of as models of $\ACF_R$ with named values for the variables $\mbf y$. By the previous case, under $\ACF_{R[\mbf y]}$ it is equivalent to a quantifier-free sentence, i.e.\ under $\ACF_R$ it is equivalent to a quantifier-free formula with free variables $\mbf y$. This concludes the proof.
\end{proof}
\end{cor}

\subsection*{Topological points}

To complete the picture, we will now show how our correspondence between models of $\ACF_R$ and geometric points descends to a bijection between elementary equivalence classes of models and topological points. Let $S_0(\ACF_R)$ denote the Stone space of $\ACF_R$, with its usual topology. There is a canonical map $p:S_0(\ACF_R)\rightarrow|\Spec(R)|$, taking a model $x:R\rightarrow k$ of $\ACF_R$ to the prime $\ker(x)$ of $R$. This is easily seen to be well-defined and continuous with respect to the constructible topology (e.g.\ the preimage of $D(f)$ is just $S_0(\ACF_R\cup\{f\neq0\})$, which is clopen in the topology on $S_0(\ACF_R)$).

We demonstrate that $p$ is a homeomorphism. That it is surjective is clear, since every topological point of $\Spec(R)$ arises from a geometric point. To prove injectivity, just note that if $x$ and $y$ are models of $\ACF_R$ such that $p(x)=p(y)$, then in particular $x$ and $y$ satisfy the same quantifier-free sentences in the language of $\ACF_R$. Quantifier elimination (corollary \ref{cor:QE_for_ACF_R}) tells us that they are elementarily equivalent, and so represent the same point of $S_0(\ACF_R)$. Thus $p$ is a continuous bijection from a compact space to a Hausdorff space, and hence is a homeomorphism.

\end{document}